\documentstyle[12pt]{article}
\begin{document}
\title{Zeros of polynomials over Cayley-Dickson algebras.}
\author{S.V. Ludkovsky.}
\date{09.06.2004}
\maketitle

\begin{abstract}
In this paper it is investigated the problem which polynomials
over Cayley-Dickson algebras (of quaternions, octonions, etc.) 
have zeros expressible with the help of roots and which
polynomials have decompositions as products of linear terms
$(z-z_{0,j})$, where $z_{0,j}$ are zeros.
\end{abstract}

\section{Introduction}
Functions of real variables with values in Clifford algebras
were investigated, for example, in \cite{brdeso}.
In this article we continue our investigations of functions
of variables belonging to noncommutative superalgebras
\cite{luoyst} considering here also functions of octonion variables
and also of more general Cayley-Dickson algebras ${\cal A}_v$,
$v\ge 4$, containing the octonion algebra as the proper subalgebra.
The Cayley-Dickson algebras of larger dimensions are not already
alternative and proceedings for them are heavier.
Nevertheless, using their power-associativity and distributibity
it is possible to investigate their polynomials.
\par Dirac had used biquaternions (complexified quaternions)
in his investigations of quantum mechanics. The Dirac
equations $D_zf_1 - D_{\tilde z}f_2=m(f_1+f_2)$ and
$D_{\tilde z}f_1 + D_zf_2=m(f_1-f_2)$ on the space
of right superlinearly $(z,{\tilde z})$-superdifferentiable
functions of quaternion variable
can be extended on the space of $(z,{\tilde z})$-superdifferentiable
functions $f_1(z,{\tilde z})$ and $f_2(z,{\tilde z})$
(see \cite{luoyst,luoyst2}) gives
evident physical interpretation of a solution $(f_1,f_2)$, $v\ge 2$, as
spinors, where $m$ is a mass of an elementary particle.
We extend the operator $(D_z^2+D_{\tilde z}^2)$ from
the space of right superlinearly $(z,{\tilde z})$-superdifferentiable
functions on the space of $(z,{\tilde z})$-superdifferentiable functions $f$, 
hence we get the Klein-Gordon equation $(D_z^2+D_{\tilde z}^2)f=2m^2f$
in the particular case of quaternions and ${\cal A}_v$-algebras, $v\ge 3$.
\par While development their theory Yang and Mills known in theoretical
and mathematical physics had actively worked with quaternions, but
they have felt lack of the available theory of quaternion functions
existing in their time. Yang also have expressed the idea, that
possibly in quantum field theory it is worthwhile to use quaternion time
(see page 198 \cite{guetze}). It is known also the
use of complex time through the Wick rotation in quantum mechanics
for getting solutions
of problems, where the imaginary time is used for interpretations
of probabilities of tunneling under energy barriers (walls).
Using the special unitary group
embedded into the quaternion skew field $\bf H$ it makes equivalent
under isomoprphism with $SO(3)$ all spatial axes. On the other hand,
the major instrument for measurement is the spectrum.
When there are deep energy wells or high energy walls, then
it makes obstacles for penetrating electromagnetic waves and radiation,
that is well known also in astronomy, where black holes are actively
studied (see page 199 and \S 3.b \cite{guetze} and references therein).
W. Hamilton in his lectures on quaternions also tackled a question of
events in $\bf H$ and had thought about use of quaternions in astronomy and
celestial mechanics (see \cite{hamilt,rothe} and references therein).
Therefore, in general to compare the sequence of events it may be
necessary in definite situations to have the same dimensional time space
as the coordinate space. On the other hand, spatial isotropy at least local
in definite domains makes from each axis under rotations and dilatations
$SU(2)\times {\bf R}$ isomorphic with $\bf H$. Therefeore, it appears
that in definite situations it would be sufficent to use ${\bf H}^4$
instead of the Minkowski space-time ${\bf R}^{1,3}$,
where ${\bf R}^{1,3}$ has the embedding into ${\bf H}$.
Since $\bf H$ as the $\bf R$-linear space is isomorphic with
$\bf R^4$, then there exists the embedding $\zeta $ of $\bf R^{1,3}$
into $\bf H$ such that $\zeta (x_1,x_i,x_j,x_k)=
x_1+x_ii+x_jj+x_kk$, where the $\bf R^{1,3}$-norm is given by
the equation $|x|_{1,3}=(x^2+{\tilde x}^2)/2=Re (x^2)=
x_1^2-x_i^2-x_j^2-x_k^2$ and the $\bf R^{1,3}$ scalar product
is given by the equality
$(x,y)_{1,3}:=(xy+{\tilde y} {\tilde x})/2=Re (xy)$
$=x_1y_1-x_iy_i-x_jy_j-x_ky_k$, where $x=x_1+x_ii+x_jj+x_kk$,
$x_1,...,x_k\in \bf R$. This also can be used
for embeddings of hyperbolic manifolds into quaternion manifolds.
Then ${\bf H}^4$ can be embedded into the algebra ${\cal A}_4$
of sedenions. It is also natural for describing systems with spin,
isospin, flavor, color and their interactions.
The enlargement of the space-time
also is dictated in some situations by symmetry properties of
differential equations or a set of operators describing a system.
For example, special unitary groups $SU(n)$, for $n=3,...,11$, etc.,
exceptional Lie groups,
are actively used in theory of elementary particles \cite{guetze}, but these
groups can be embedded into the corresponding
Cayley-Dickson algebra ${\cal A}_v$ \cite{baez}.
Indeed, $U(m)\subset GL(n,{\bf C})\subset {\bf C}^{n^2}$ while
${\bf C}^m$ has the embedding into ${\cal A}_r$, where
$m=2^{r-1}$, such that ${\bf C}^m \ni (x^1+iy^1,...,
x^m+iy^m)=:\xi \mapsto z:= (x^1+i_1y^1)+i_2(x^2+i_2^*i_3y^2)+...+
i_{2m-2}(x^m+i_{2m-2}^*i_{2m-1}y^m)\in {\cal A}_r$, since
$(i_l^*i_k)^2=-1$ for each $l\ne k\ge 1$, where $ \{ i_0,i_1,...,
i_{2m-1} \} $ is the basis of generators of ${\cal A}_r$, $i_0=1$,
$i_k^2=-1$, $i_0i_k=i_ki_0$, $i_li_k=-i_ki_l$ for each $k\ne l\ge 1$,
$i=(-1)^{1/2}$, $z^*=x^1-i_1y^1-i_2x^2-i_3y^2-...-i_{2m-2}x^m-
i_{2m-1}y^m$, the norm $(zz^*)^{1/2}=:|z|=(\sum_{k=1}^m
|x^k+iy^k|^2)^{1/2}=:|\xi |$ satisfies the parallelogramm
identity and induces the scalar product.
\par This paper continuous investigations of function theory
over Cayley-Dickson algebras \cite{luoyst,luoyst2}.
Cayley-Dickson algebras ${\cal A}_v$ over the field of real numbers
coincide with the field $\bf C$ of complex numbers for $v=1$,
with the skew field of quaternions, when $v=2$,
with the division nonassociative noncommutative algebra $\bf K$ of octonions
for $v=3$, for each $v\ge 4$ they are nonassociative and not
division algebras. The algebra ${\cal A}_{v+1}$ is obtained from
${\cal A}_v$ with the help of the doubling procedure.
Considering the completion of the limit of the direct system
of Cayley-Dickson algebras $ \{ {\cal A}_v: v\in {\bf N} \} $
relative to the $l_2$ norm, we get the Cayley-Dickson algebra
${\cal A}_{\infty }$ with infinite number of generators (see
\cite{luoyst2}). Denote by ${\cal I}_v$ the set of all
$z\in {\cal A}_v$ such that $Re (z)=0$, where $Re (z):= (z + z^*)/2$.
In the paper \cite{huso} there were written formulas for zeros
of square polynomials of particular type over quaternions.
\par In this paper it is investigated the problem which polynomials
over Cayley-Dickson algebras ${\cal A}_v$ with $v\ge 2$
have zeros expressible with the help of roots and which
polynomials have decompositions as products of linear terms
$(z-z_{0,j})$, where $z_{0,j}$ are zeros.
\par The results of this paper can serve for subsequent investigations
of special functions of Cayley-Dickson algebra variables,
noncommutative sheaf theory, manifolds of noncommutative geometry
over Cayley-Dickson algebras, their groups of loops and diffeomorphisms
(see also \cite{lulgcm,lulsqm,lustptg,lupm,oystaey}).
\section {Zeros of polynomials.}
\par {\bf 1. Remark.} Let $P_n(z)$ be a polynomial over the
Cayley-Dickson algebra ${\cal A}_v$, $v\ge 2$, such that
$$P_n(z)=\sum_{0\le |k|\le n} (a_k,z^k),$$
where $a_k=(a_{k_1},...,a_{k_m})$, $m=m(k)\in \bf N$,
$a_{k_1},...,a_{k_m}\in {\cal A}_v$, $z\in {\cal A}_v$,
$(a_k,z^k):= \{ a_{k_1}z^{k_1}...a_{k_m}z^{k_m} \} _{q(2m)}$,
$k=(k_1,...,k_m)$, $0\le k_l\in \bf N$ for each $l=1,...,m$,
$|k|=k_1+...+k_m$, $q(m)$ is the vector indicating on the order
of subsequent multiplications in  $ \{ * \} $ \cite{luoyst2}.
A number $z_0\in {\cal A}_v$ is called as usually
a zero of $P_n$ if $P_n(z_0)=0$. Denote by $deg (P_n):=
\max \{ |k|:$ $\mbox{there exist}$ $a_{k_1}\ne 0,...,a_{k_m}\ne
0, k=(k_1,...,k_m) \} =: n$ the degree of $P_n$.
\par {\bf 2. Definition.} Let $M$ be a set such that $M=\bigcup_jU_j$,
$M$ is a Hausdorff topological space, each $U_j$ is open in $M$,
$\phi _j: U_j\to \phi _j(U_j)\subset {\cal A}_v^m$ are homeomorphisms,
$\phi _j(U_j)$ is open in ${\cal A}_v$ for each $j$, if $U_i\cap U_j\ne
\emptyset $, the transition mapping $\phi _i\circ \phi _j^{-1}$
is ${\cal A}_v$-holomorphic on its domain. Then $M$ is called
the ${\cal A}_v$-holomorphic manifold and $m\in \bf N$ is its dimension
over ${\cal A}_v$.
\par {\bf 3. Definition.} Let $G$ be a group (may be without
associativity axiom) such that $G$ has also the structure of an
${\cal A}_v$-holomorphic manifold and the mapping
$G^2\ni (g,h)\mapsto gh^{-1}\in G$ is continuous (${\cal A}_v$-holomorphic),
then $G$ we call the Cayley-Dickson-Lie group (${\cal A}_v$-holomorphic).
We also suppose that $gh=0$ is not excluded for some $g, h\in G$,
when $v\ge 4$, where $0G=G0= \{ 0 \} $, $G$ is an additional zero element,
while $G$ is considered as the multiplicative group with
the unit element $1$.
\par {\bf 4. Definition.} Let $P_n$ be a polynomial over ${\cal A}_v$,
$v\ge 2$. Denote by ${\cal G}_l(P_n,0):= \{
g\in {\cal A}_v^{\bullet }:$ $\mbox{there exist}$ $z_0$ $\mbox{and}$
$z_1$ $\mbox{with}$ $P_n(z_0)=0$ $\mbox{and}$ $P_n(z_1)=0$
$\mbox{such that}$ $gz_0=z_1$ $\mbox{and}$ $g^{-1}z_1=z_0 \} ,$
${\cal G}_l(P_n) := \{ g\in {\cal A}_v^{\bullet }:$ $P_n(gz)=P_n(z)
\forall z \in {\cal A}_v \} $, analogously ${\cal G}_r(P_n,0)$
and ${\cal G}_r(P_n)$ are defined with respect to the right action.
\par {\bf 5. Lemma.} {\it Let $P_n$, ${\cal G}_l(P_n,0)$
and ${\cal G}_l(P_n)$ be as in Definition 4, then ${\cal G}_l(P_n,0)$
and ${\cal G}_l(P_n)$ are subgroups of the multiplicative
Cayley-Dickson-Lie group ${\cal A}_v^{\bullet }$, where
$2\le v\le \infty $.}
\par {\bf Proof.} Suppose that $0\in {\cal G}_l(P_n)\times {\cal G}_l(P_n)$,
that is, there exist $g, h\in {\cal G}_l(P_n)$ such that $gh=0\in
{\cal A}_v$. Consider $z\in {\cal A}_v$ such that $hz\ne 0$ and
$h^{-1}(hz)=z$, which is possible due to the embedding
${\bf K}\hookrightarrow {\cal A}_v$. Then
also $h(h^{-1}z)=z$. Therefore, $P_n(g(h(h^{-1}z))=
P_n(gz)=P_n(z)$. Take any $z\in \bf R$, then
$P_n(g(hz))=P_n((gh)z)=P_n(0)$, consequently, $P_n(z)=P_n(0)$
for each $z\in \bf R$. This means that $deg P_n(z)=0$, since
$P_n(0)=a_0$. Therefore, ${\cal G}_l(P_n)={\cal A}_v^{\bullet }$.
Evidently, ${\cal A}_v^{\bullet }$ is the Cayley-Dickson-Lie
group the atlas of which can be chosen consisting of one chart.
Mention, that this case is impossible for $v=2$ and $v=3$
and nontrivial $P_n$, since $\bf H$ and $\bf K$ are division algebras.
\par Suppose now that $0\notin {\cal G}_l(P_n)\times {\cal G}_l(P_n)$,
then ${\cal G}_l(P_n)$ has not divisors of zero, consequently,
${\cal G}_l(P_n)$ can be embedded into $\bf K$.
For $g, h\in {\cal G}_l(P_n)$ choose $z\in \Upsilon _{M_g,M_h}
\hookrightarrow {\bf H}$, where $\Upsilon _{M_g,M_h}$ is the minimal
subalgebra of ${\cal A}_v$ containing $M_g$ and $M_h$ and hence
$g$ and $h$ also, $g=\rho _g\exp (M_g)$ is the polar decomposition
with $\rho _g\ge 0$ and $M_g\in {\cal I}_v$,
then $P_n(g(hz))=P_n((gh)z)=P_n(z)$ for each $z\in \Upsilon _{g,h}$,
since $\bf H$ is associative. If $|g|\ne 1$ and $deg (P_n)\ge 1$, then from
$P_n(g^kz)=P_n(z)$ for each $z\in \Upsilon _{M_g,M_{\xi }}={\bf H}
\hookrightarrow {\cal A}_v$, where $|z|\ne 0$ in particular can be taken,
$k\in \bf Z$ is arbitrary it follows, that $\lim_{k\to \infty }
P_n(g^kz)=0$ for $|g|>1$ and $\lim_{k\to - \infty }
P_n(g^kz)=0$ for $|g|<1$, consequently, $P_n(z)=0$ for each
$z\in {\cal A}_v$, that contradicts our assumption. It remains
the case $|g|=1$ for each $g\in {\cal G}_l(P_n)$, then
$g=\exp (M_g)$ with $M_g\in {\cal I}_v$ due to the polar decomposition
of Cayley-Dickson numbers. Therefore, there exists the embedding
$\Upsilon _{M_g,M_h,M_z}\hookrightarrow {\bf K}\hookrightarrow
{\cal A}_v$ for corresponding $M_g,M_h,M_z$ (see also \cite{kansol})
for each $z=\rho _z\exp (M_z)\in {\cal A}_v$ with
$\rho _z\ge 0$ and $M_z\in {\cal I}_v$, since $\Upsilon _{M_g,M_h}
\hookrightarrow \bf H$. But $g{\bf K}=\bf K$ for each
$g\ne 0$, $g\in \bf K$, consequently, $P_n(g(hz))=P_n(hz)=P_n(z)$
$=P_n((gh)(h^{-1}\eta ))=P_n(h^{-1}\eta )=P_n(\eta )$ for each
$z\in \Upsilon _{M_g,M_h,M_z}\hookrightarrow \bf K$
due to the alternative property of $\bf K$, where $\eta :=hz$.
If $hz=0$, then $h^{-1}z=0$ and $P_n((gh)(h^{-1}z))=P_n(0)=P_n(h^{-1}z)$.
In general, $z=z_1+z_2$, where $z_1\in \Upsilon _{M_g,M_h,M_{z_1}}
\hookrightarrow \bf K$, $hz_2=0$, then $P_n((gh)(h^{-1}z))=
P_n((gh)(h^{-1}z_1))=P_n(h^{-1}z_1)=P_n(h^{-1}z)=P_n(z)$.
Then $P_n(g(g^{-1}z))=P_n(g^{-1}z)=P_n(z)$ for each
$z\in {\bf K}_g$, where ${\bf K}_g$ is any copy of $\bf K$ containing
$g\in {\cal A}_v$, $gg^{-1}=1$, $P_n(1z)=P_n(z)$ for all $z\in {\cal A}_v$.
Therefore, $g^{-1}\in {\cal G}_l(P_n)$ for each $g\in {\cal G}_l(P_n)$
and $1$ is the unit element of ${\cal G}_l(P_n)$,
hence ${\cal G}_l(P_n)$ is the group by Definition 3.
\par Consider now ${\cal G}_l(P_n,0),$ then $P_n(g(hz_0))=
P_n(hz_0)=P_n(z_0)=0$, $P_n(g^{-1}z_1)=P_n(z_0)=0$
for some zeros $z_0$ and $z_1$ ($P_n(z_1)=0$ and $P_n(z_0)=0$)
of $P_n$ with $gz_0=z_1$ and $g^{-1}z_1=z_0$. Thus
$g^{-1}\in {\cal G}_l(P_n,0)$ for each $g\in {\cal G}_l(P_n,0)$
and $gh\in {\cal G}_l(P_n,0)$ for each $g, h\in {\cal G}_l(P_n,0)$.
\par {\bf 6. Corollary.} {\it If $v=2$, then ${\cal G}_l(P_n)$ and
${\cal G}_l(P_n,0)$ are associative, if  $v=3$, then they are alteranative,
in both cases $v=2$ and $v=3$ they have not divisors of zero.}
\par {\bf 7. Corollary.} {\it ${\cal G}_l(P_n)\subset {\cal G}_l(P_n,0)$,
${\cal G}_l(P_nP_m)\supset {\cal G}_l(P_n)\cap {\cal G}_l(P_m)$
for each $P_n$ and $P_m$ over ${\cal A}_v$.}
\par {\bf 8. Definition.} We say that a zero $z_0$ of $P_n$
is expressed in roots, if there exist $n_j\ge 2$, polynomials
$f_p$ and a polynomial $\xi $ such that
$$z_0=\xi (f_{i(1)}^{1/n_{i(1)}}\circ f_{i(2)}^{1/n_{i(2)}}\circ ...
\circ f_{i(k)}^{1/n(i(k))},...,f_{i(w)}^{1/n_{i(w)}}\circ
f_{i(w+1)}^{1/n_{i(w+1)}} \circ ... \circ f_{i(t)}^{1/n(i(t))}),$$
where $i(1),...,i(t)\in \{ 1,...,p \} $, $f_i$ are polynomials
of $(a_k: k)$ or also of others polynomials or their roots,
all roots are certainly with $1\le n_i\le n$.
\par {\bf 9. Lemma.} {\it Let $M, N\in {\cal I}_v$,
$N=N_1+N_2$, where $M\perp N_2$
relative to the scalar product $(x,y) := Re (xy^*)$, $v\ge 2$,
$N_1=\beta M$, where $\beta \in \bf R$, then
the commutator of $e^M$ and $e^N$ is:
\par $[e^M,e^N]=2[(\sin |M|)/|M|] [(\sin |N|)/|N|]MN_2$.}
\par {\bf Proof.} Since $Re (MN^*)=0$ and $N^*=-N$, then
$Re (MN)=0$. We have $e^Ne^M=(\cos |M|)e^N+[(\sin |M|)/|M|]M
e^{N_1-N_2}$ and $e^Me^N=(\cos |M|)e^N+[(\sin |M|)/|M|]Me^N$
(see Formula $3.3$ in \cite{luoyst2}).
Therefore, $e^Me^N-e^Ne^M=[(\sin |M|)/|M|] [e^N-e^{N_1-N_2}]$.
On the other hand, $|N|^2=-N_1^2-N_2^2=|N_1-N_2|^2$
and $e^{N_1-N_2}=\cos |N| + [(\sin |N|)/|N|] (N_1-N_2)$,
consequently, $e^N - e^{N_1-N_2}=2 [(\sin |N|)/|N|] N_2$.
\par {\bf 10. Corollary.} {\it If $M, N\in {\cal I}_v$, $v\ge 2$,
then $[e^M,e^N]=0$ if and only if
$|N|=\pi l$ or $|M|=\pi l$ or $M\parallel N$, where $0\le l\in \bf Z$.}
\par {\bf Proof.} The condition $M\parallel N$ is equivalent to
$MN_2=0$.
\par {\bf 11. Lemma.} {\it Let $M, N \in {\cal I}_v$ and $M\perp N$
relative to the scalar product $(x,y) := Re (xy^*)$, $v\ge 2$,
then the commutator of $e^M$ and $e^N$ is:
\par $[e^M,e^N]=2[(\sin |M|)/|M|][(\sin |N|)/|N|]MN$.}
\par {\bf Proof.} Since $e^Ne^M=(\cos |M|) e^N +[(\sin |M|)/|M|]M
e^{-N}$, consequently,
$[e^M,e^N]=e^Me^N-e^Ne^M=[(\sin |M|)/|M|] M (e^N-e^{-N})$,
where $e^N-e^{-N}=2[(\sin |N|)/|N|]N$.
\par {\bf 12. Lemma.} {\it Let $M, N\in {\cal I}_v$,
$N=N_1+N_2$, where $M\perp N_2$, $v\ge 2$,
$N_1=\beta M$, where $\beta \in \bf R$, then
the anticommutator of $e^M$ and $e^N$ is:
\par $\{ e^M, e^N \} =2(\cos |M|)e^N + 2 [(\sin |M|)/|M|] M
(\cos |N|+ [(\sin |N|)/|N|]N_1)$.}
\par {\bf Proof.} As follows from \S 9: \\
$\{ e^M,e^N \} := e^Me^N+e^Ne^M=
2(\cos |M|) e^N + [(\sin |M|)/|M|]M[e^N+e^{N_1-N_2}]$,
where $e^N+e^{N_1-N_2}=2\cos |N|+2[(\sin |N|)/|N|]N_1$.
\par {\bf 13. Corollary.} {\it If $M, N\in {\cal I}_v$, 
then $ \{ e^M, e^N \} =0$ if and only if either $\{ |N_2|=\pi /2+\pi l$ and
$N_1=0$ and $|M|=\pi /2 +\pi k \} $ or $ \{ |M_2|=\pi /2+\pi k$ and
$M_1=0$ and $|N|=\pi /2+\pi l \} $, where $0\le l, k \in \bf Z$,
$N_1=\beta M$ and $N_2\perp M$, $M_1=\alpha N$ and $M_2\perp N$,
$N=N_1+N_2$, $M=M_1+M_2$, $\alpha , \beta \in \bf R$.}
\par {\bf 14. Lemma.} {\it Let $K\in {\cal I}_v$, $\infty \ge v\ge 2$,
then there exists $N\in {\cal I}_v$ such that $e^N(e^Ke^N) \in \bf R$.
Moreover, the manifold of such $N$ has the real codimension $2$
for $e^K\notin \bf R$ and $1$ for $e^K\in \bf R$.}
\par {\bf Proof.} For $|K|=0\quad (mod \pi )$ this is trivial.
Consider $|K|\ne 0 \quad (mod \pi )$.
At first $e^Ke^N=e^N\cos |K| + e^{N_1-N_2}
[(\sin |K|)/|K|]K$, where $N=N_1+N_2$, $N_1=\beta K$, $N_2\perp K$,
$\beta \in \bf R$. Then \\
$e^N(e^Ke^N)=e^{2N} \cos |K| + e^N(e^{N_1-N_2} [(\sin |K|)/ |K|]K)$. \\
On the other hand, $e^N= \cos |N| + [(\sin |N|)/|N|]N$, 
\par $e^{N_1-N_2}= \cos |N| + [(\sin |N|)/|N|](N_1-N_2)$, \\
$e^N e^{N_1-N_2}= \cos ^2|N| +
[(\sin ^2 |N|)/|N|^2] (N_1^2-N_2^2+2N_2N_1)
+ (\cos |N|) [(\sin |N|)/|N|] 2N_1$, \\
since $N_1\perp N_2$ and $N_1N_2=-N_2N_1$, where 
$N_1N_2\in {\cal I}_v$, since $N_1\perp N_2\in {\cal I}_v$.
Since $Re (KN_2)=0$, then
\par $Im (e^{-N}(e^Ke^N))=0$ if and only if
\par $[(\sin |2N|)/|N|]N \cos |K| + $ \\
$(\cos ^2|N|+[(\sin ^2 |N|)/|N|^2] (N_1^2-N_2^2+2N_2N_1)
[(\sin |K|)/ |K|]K=0$. \\
Evidently, the manifold $\Psi (K)$ of $N\in {\cal I}_v$ satisfying these
conditions has $codim_{\bf R}\Psi (K)=2$ for each $v\ge 2$.
\par {\bf 15. Lemma.} {\it Let $n\ge 2$ be a natural number,
$r>0$, then the manifold $\Phi _n(r)$ of all $\zeta \in {\cal A}_v$
such that $\zeta ^n=r$ has $codim_{\bf R} \Phi _n(r)=2$ for each
$2\le v\le \infty $.}
\par {\bf Proof.} Represent $\zeta \in \Phi _n(r)$ as the product
$\zeta =\rho z$, where $\rho >0$ and $z\in \Phi _n(1)$,
$\rho ^n=r$. Therefore, it is sufficient to describe $\Phi _n(1)$.
Consider $M\in {\cal I}_v$ such that $(e^M)^n=e^{nM}=1$.
In view of Formulas $(3.2,3)$ \cite{luoyst2} it is the case
if and only if $|M|=2\pi k$, where $0\le k\in \bf Z$.
This defines manifold $\Phi _n(1)$ of codimension $2$ in ${\cal A}_v$.
\par {\bf 16. Lemma.} {\it Let $\zeta \in {\cal A}_v$, $2\le v\le \infty $,
$n\ge 2$, $z\ne 0$, then the manifold $\Phi _n(\zeta )$ of all
$z \in {\cal A}_v$ such that $z^n=\zeta $ has
$codim_{\bf R}\Phi _n(\zeta )=2$ and $\Phi _n(\zeta )$ is locally connected.}
\par {\bf Proof.} In the polar form $\zeta =re^K$, where
$r>0$, $K\in {\cal I}_v$. In view of Lemma 15 it remains to consider
$\zeta \in {\cal A}_v\setminus [0,\infty )$. In accordance with Lemma
14 there exists the $2$-codimensional manifold of $N\in {\cal I}_v$
such that $e^{-N}(e^Ke^{-N}) \in \bf R$. In particular, take
$N=K/2$. Represent each
$z\in \Phi _n(\zeta )$ in the form $z=\rho e^N(e^Me^N)$,
where $\rho >0$ and $\rho ^n=r$. For this verify that
each $0\ne z\in {\cal A}_v$ can be written in such form, when
$N$ is given. For this consider the minimal subalgebra
$\Upsilon _{L,N}$ of ${\cal A}_v$ containing $L$ and $N$,
where $z=\rho e^L$, $L\in {\cal I}_v$. For the corresponding $N$
such that $N\perp L$ and $|NL|=|N| |L|$
we have the embedding $\Upsilon _{L,N}\hookrightarrow {\bf H}$.
Therefore, $e^M=\rho ^{-1} (e^{-N}z)e^{-N}$ from which the
existence of $M\in \Upsilon _{L,N}$ follows, where $\rho :=|z|$,
since $|e^N|=1$ and hence $|(e^{-N}z)e^{-N}|=\rho $.
Therefore, we have the equation $(e^N(e^Me^N))^n=e^K$.
Since $codim_{\bf R}\Phi _n(e^{-N}(e^Ke^{-N})) = 2$
and for each $z\in \Phi _n(e^{-N}(e^Ke^{-N}))$ there exists
$M\in {\cal I}_v$ such that $e^N (e^Me^N) =z$.
Then $codim_{\bf R}\Phi _n(\zeta )=2$.
From the formulas above it follows, that $\Phi _n(\zeta )$
is locally connected.
\par {\bf 17. Theorem.} {\it Let $P_n$ be a polynomial, $n\ge 2$,
$z_0$ be its zero expressed in roots, then the connected component
$S(z_0)$ of $z_0$ in the set of all zeros ${\cal O}(P_n):=
\{ z\in {\cal A}_v:$ $P_n(z)=0 \} $ is a positive dimensional
real algebraic submanifold of ${\cal O}(P_n)$ embedded into
${\cal A}_v$, $2\le v\in \bf N$.
The group ${\cal G}_l(P_n,0)$ is the locally connected
unique up to isomorphism subgroup of ${\cal A}_v^{\bullet }$ and
$codim_{\bf R} {\cal G}_l(P_n,0)\le 2$.}
\par {\bf Proof.} In view of Lemma 16 $codim_{\bf R}\Phi _n(\zeta )=
2$ for each $0\ne \zeta \in {\cal A}_v$, $n\ge 2$. The application
of $k$-roots to $\Phi _n(\zeta )$ with $k\ge 2$ gives
$codim_{\bf R}\Phi _k(z)=2$ for each $z\in \Phi _n(\zeta )$,
consequently, $codim_{\bf R} \Phi _k(\Phi _n(\zeta ))\le 2$.
Therefore, for each zero $z_0\in {\cal O}(P_n)$ we have
$codim_{\bf R}S(z_0)\le 2$ and $S(z_0)$ is connected.
Consider a local subgroup
${\cal G}_l(P_n,0;\zeta _0)$ of ${\cal A}_v^{\bullet }$ elements of
which relate all corresponding pairs $z_0, z_1 \in S(\zeta _0)$.
In view of Lemma 16 for $dim_{\bf R}(S(\zeta _0))>0$ and
$dim_{\bf R}(S(\zeta _1))>0$ for each $\zeta _0\ne \zeta _1\in
{\cal A}_v$ there exists the embedding ${\cal G}_l(P_n,0;\zeta _i)
\hookrightarrow {\cal G}_l(P_n,0;\zeta _j)$ defined up to
group isomorphisms for either
$(i=1,j=2)$ or $(i=2,j=1)$, since $e^{-N}((\Phi _n(z))e^{-N})=
\Phi _n(|z|)$ for $z\in {\cal A}_v$ and the corresponding
$N\in {\cal I}_v$. In view of formulas for $\Phi _n(z)$
these local subgroups generate the group ${\cal G}_l(P_n,0)$
after addition transition elements between corresponding
connected components arising from the same algebraic expressions
in radicals. The group ${\cal G}_l(P_n,0)$ exists in
accordance with Lemma 5. Evidenly, ${\cal G}_l(P_n,0)$ is locally
connected and $codim_{\bf R}{\cal G}_l(P_n,0)\le 2$.
\par {\bf 18. Theorem.} {\it Let $G$ be a locally connected subgroup
of the unit ball $B:=B({\cal A}_v^{\bullet },0,1)$ with centre $0$,
$2\le v\in \bf N$,
then there exists a nonnegative right-invariant $\sigma $-additive
measure $\mu $ on the Borel algebra ${\cal B}(G)$ of $G$.}
\par {\bf Proof.} Since $|z \zeta | \le |z| |\zeta |$ for each
$z, \zeta \in {\cal A}_v$, then $B$ is the multiplicative
group in accordance with Definition 2. There exists the retraction
$f$ of $B$ on $G$ such that $f(z)=z$ for each $z\in G$ and $f(B)=G$.
Suppose that $\nu $ is a right-invariant measure
$\nu : {\cal B}(B) \to [0,1]\subset \bf R$.
Then it induces the measure $\psi (A):=\nu (f^{-1}(A))$ for
each $A\in {\cal B}(G)$. It is right-quasi-invariant:
$\psi ((dz) g^{-1})=\rho (g,z) \psi (dz)$ for each $g, z\in G$,
where $\rho (g,z) \in L^1(G,{\cal B}(G),\mu ,{\bf R})$.
Put $\mu (dz):=\rho ^{-1}(z,1) \psi (d\xi )$, $d\xi :=(dz)z^{-1}$,
then $\mu ((dz)g^{-1}) = \mu (dz)$, since $\mu (
\{ hg:$ $h\in A, hg=0 \} )=0$ for each $A\in {\cal B}(G)$.
It remains to construct $\nu $. For this take
the Hermitian metric $z\zeta ^* =: g(z,\zeta )$ in ${\cal A}_v$.
Relative to it $e^M$ gives the unitary operator for the quadratic
form $g(z,z)=g(ze^M,ze^M)$ for $\lambda $-almost every $z\in {\cal A}_v$,
preserving the measure on ${\cal A}_v$ induced by the
Lebesgue measure $\lambda $ on the real shadow ${\bf R}^{2^v}$.
Then $\lambda ((dz)g) = |g|\lambda (dz)$, so put
$\nu (dz) := |z|^{-1}\lambda (dz)$ on ${\cal B}({\cal A}_v^{\bullet })$,
hence $\nu ((dz)g^{-1})=\nu (dz)$, since
$\nu (S)=0$ for each submanifold $S$ of ${\cal A}_v$ with
$codim_{\bf R}S\ge 1$ and  $\nu ( \{ hg:$  $h\in A, hg=0 \} )=0$.
\par {\bf 19. Definition.} Let $G$ be a subgroup of $B$ as
above. If $ \{ e^{tM}:$ $t \in [0, 2\pi ] \} $ is a one-parameter
subgroup of $B$, then take $t\in [0,\pi /(2n)]$. With this restriction
consider a subset $H_n$ containing intersections with $G$ of such
local subgroups instead of the entire one-parameter subgroup
such that $H_n^{4n}=G$, where $H_n^2=H_nH_n$,...,$H_n^{k+1}=
H_n^kH_n$ for each $k\in \bf N$, $AB:=\{ ab: a\in A, b\in B \} $
for two subgroups $A$ and $B$ of a given group $G$.
If $P_n$ is a polynomial on ${\cal A}_v$,
then put $P_n^G(z) := \int_{H_n} P_n(gz)\mu (dg)$, where $\mu $
is the right-invariant measure on $G$.
\par {\bf 20. Theorem.} {\it Let $P_n$ be a polynomial on
${\cal A}_v$, $2\le v\in \bf N$. Then $P_n^G(z)=P_n^G(1)$ for each $z\in G$
(see \S 19).
\par If $v=2$, then ${\cal G}_l(P_n^G)\supset G$ and
$(P_n^G)^G=cP_n^G$, where $c=const >0$.
\par If $G\subset {\cal G}_l(P_n,0)$, then ${\cal G}_l(P_n^G,0)
\supset {\cal G}_l(P_n,0)$.}
\par {\bf Proof.} There are integral equalities: \\
$P_n^G(z):=\int_{H_n} P_n(gz)\mu (dg)
=\int_{H_n} P_n(y)\mu ((dy)z^{-1})=\int_{H_n} P_n(y)\mu (dy)=P_n^G(1)$ \\
for each $z\in H_n$. Iterating this equations we get the statement
of this theorem, since $H_n^{4n}=G$.
Verify that this construction gives nontrivial in genereal
polynomial. Indeed $\int_{H_n} g^k \mu (dg)\ne 0$
for each $0\le k\le n$ by the definition of $H_n$.
\par In the case of $\bf H$, $v=2$, we have \\
$P_n^G(xz)= \int_{H_n}P_n(g(xz))\mu (dg)=
\int_{H_n}P_n(yz)\mu ((dy)x^{-1})$ \\
$=\int_{H_n}P_n(yz)\mu (dy)=P_n^G(z)$ for each $x\in H_n$, \\
iterating this equation $1\le k\le 4n$ times,
we get the statement for each $x\in G$. In view of the Fubini
theorem \\
$(P_n^G)^G(z)=\int_{H_n} \int_{H_n} P_n(g_1g_2z)\mu (dg_1)\mu (dg_2)$ \\
$=\int_{H_n}(\int_{H_n}P_n(hz)\mu ((dh)g_2^{-1}))\mu (dg_2)$ \\
$= \int_{H_n} P_n(hz) \mu (dh)=c P_n^G(z)$, \\
since $\mu $ is right-invariant probability measure,
\par $\int_{H_n} \mu (Sg^{-1}) \mu (dg)=\mu (S) \mu (H_n)$ \\
for each $S\in {\cal B}(G)$ and $\mu (H_n)=:c>0$,
$P_n (z) = \sum_k P_{n,k}(z) i_{2k-1}$, where $P_{n,k}(z)$ is with values
in ${\bf R}\oplus {\bf R}i_{2k} i^*_{2k-1}$ isomorphic with $\bf C$,
$\{ 1, i_1,...,i_{2^v-1} \} $ are generators of ${\cal A}_v$.
\par If $G\subset {\cal G}_l(P_n,0)$ and $z_0\in {\cal O}(P_n)$,
$g\in {\cal G}_l(P_n,0)$, then
\par $P_n^G(gz_0)=\int_{H_n}P_n(g(xz_0))\mu (dx)=0$. 
\par It remains to verify that $P^G(z)$ is also polynomial,
if $P$ is a polynomial. Mention, that if $P_k$ is
$\bf R$-homogeneous of degree $k$, then after integration
the function $P^G_k$ is also $\bf R$-homogeneous:
$P^G_k(az)=a^kP^G_n(z)$ for each $a\in \bf R$, since the integral
$\int_{H_n}P_k(gz)\mu (dg)$ is $\bf R$-linear.
Therefore, in general $P^G_n$ is the sum of $\bf R$-homogeneous terms
of degree at most $n$. The set of all step functions
$\phi (g) := \sum_{j=1}^m \chi _{A_j}(g) f(g_jz)$ is dense in the space
$L^1 (G, {\cal B}(G), \mu , {\cal A}_v)$ of ${\cal A}_v$-valued
$\mu $-measurable functions of finite norm $\| f \| _{L^1}:=
\int_G|f(z)| \mu (dz)$, where $|\eta |$ is the norm of
$\eta $ in ${\cal A}_v$, $\chi _A$ is the characteristic function
of $A\in {\cal B}(G)$, $g_j\in A_j$. Since $H_n\in {\cal B}(G)$, then
$\int_{H_n} \sum_{j=1}^m \chi _{A_j}(g) P_n(g_jz)\mu (dg)=
P_n(g_jz)\mu (A_j\cap H_n)$ is the polynomial function
of degree $n$ by the variable $z\in {\cal A}_v$. Therefore,
the limit relative to $\| * \| _{L^1}$ of the Cauchy sequence
of such step functions and their integrals by $H_n$ give a sequence
converging to the polynomial $P^G_n(z)$ of degree $n$.
This convergence is uniform by the variable $z$ in each
ball $B({\cal A}_v,0,r)$ with $0< r<\infty $, that follows from
the consideration of integrals of step functions and since
$P_n$ is characterized by the finite family of expansion
coefficients.
\par {\bf 21. Remarks.} To construct polynomials with prescribed
dimensions of $S(z_0)$, which can vary from $0$ to $2^v-1$
($2^v$ corresponds to the trivial case),
for $z_0\in {\cal O}(P_n)$ it is possible
to use projection operators on $span_{\bf R} (i_{j_1},...,i_{j_k})$
from \cite{luoyst2} and as well expressions of commutators
and anticommutators above, where $(1, i_1,..., i_{2^v-1})$
are generators of ${\cal A}_v$, $2\le v\in \bf N$. Indeed, each
projection $\pi _j$ on ${\bf R}i_j$ has the form $\pi _j(z)=
\sum_k \alpha _{j,k} (z \beta _{j,k})$, where $\alpha _{j,k}$ and
$\beta _{j,k}$ are constants in ${\cal A}_v$ independent from
$z\in {\cal A}_v$.
Therefore, in general a polynomial over ${\cal A}_v$ may have
a set of zeros different from discrete, hence the representation
of the type $a (z-z_1)...(z-z_n)$ is the exceptional case for $v\ge 2$,
moreover, for $v\ge 3$ because of nonassociativity of ${\cal A}_v$,
where $z_1,...,z_n$ are zeros.
Theorems above also show, that in general a polynomial $P_n$
over ${\cal A}_v$, $v\ge 2$, with $n\le 4$ may have not expressions
of zeros through roots apart from the complex case.
\par Recall that in accordance with the Abel theorem the general equation
over $\bf C$ of degree $n$ is solvable in radicals if and only if $n\le 4$
(see Proposition 9.8 on page 308 in \cite{hung}).
Since the general polynomial $P_n(z)$ over ${\cal A}_v$ has a restriction
on $\bf C$ embedded in ${\cal A}_v$, then the general equation
of degree $n$ over ${\cal A}_v$, $v\ge 2$, is solvable in radicals
if and only if $n\le 4$, but this causes also the restriction on dimensions
of $S(z_0)$ (see above).
\par For ${\cal A}_{\infty }$ it is possible to define as in \S 19
$P_n^G(z):= \lim_{m\to \infty } P_n^{G_m}(z_m)$, where
$G_m:=G\cap {\cal A}_m$, $z_m$ is the projection of $z\in {\cal A}_{\infty }$
into ${\cal A}_m$.
\par {\bf 22. Corollary.} {\it Let $P_n$ be a polynomial on
${\cal A}_{\infty }$. Then $P_n^G(z)=P_n^G(1)$ for each $z\in G$.
\par If $G\subset {\cal G}_l(P_n,0)$, then ${\cal G}_l(P_n^G,0)
\supset {\cal G}_l(P_n,0)$. }
\par {\bf Proof.} For each $g\in G$ and $z\in {\cal A}_{\infty }$
the sequences $g_m\in G_m$ and $z_m \in {\cal A}_m$ converge to
$g$ and $z$ respectively, while $m$ tends to the infinity.
Therefore, $P_n^G(z)=\lim_{m\to \infty }P_n^{G_m}(z_m)=
\lim_{m\to \infty } P_n^{G_m}(1)=P_n^G(1)$ for each $z\in G$.
\par For each $z\in {\cal A}_{\infty }$ we have $P_n^G(z)=\lim_{m\to \infty }
P_n^{G_m}(z_m)=\lim_{m\to \infty } P_{n,m}^{G_m} (z_m)$, where
$P_{n,m}$ is the polynomial with coefficients $(a_{k_j})_m\in {\cal A}_m$
instead of $a_{k_j}\in {\cal A}_{\infty }$.
Thus, for each $g\in {\cal G}_l(P_n^G,0)$ there exists a sequence
$g_m\in {\cal G}_l(P_{n,m}^{G_m},0)$ converging to $g$,
when $m$ tends to the infinity.
\par The following theorem establishes, that each polynomial over
Cayley-Dickson algebra has zeros.
\par {\bf 23. Theorem.} {\it Let $P(z)$ be a polynomial
on ${\cal A}_v$, $2\le v\le \infty $,
such that $P(z)=z^{n+1}+\sum_{\eta (k)=0}^n(A_k,z^k)$,
where $A_k=(a_{1,k},...,a_{m,k}),$ $a_{j,l}\in {\cal A}_v,$
$k=(k_1,...,k_m)$, $0\le k_j\in \bf Z$, $\eta (k)=k_1+...+k_m$,
$0\le m=m(k)\in \bf Z$, $m(k)\le \eta (k)+1$,
$(A_k,z^k):=\{ a_{1,k}z^{k_1}...a_{m,k}z^{k_s} \} _{q(m+\eta (k))}$,
$z^0:=1$. Then $P(z)$ has a zero in ${\cal A}_v$.}
\par {\bf Proof.} Consider at first $v<\infty $.
Suppose that $P(z)\ne 0$ for each $z\in {\cal A}_v$.
Consider a rectifiable path $\gamma _R$ in ${\cal A}_v$ such that
$\gamma _R ([0,1]) \cap {\cal A}_v = [-R,R]$ and outside $[-R,R]$:
$\quad \gamma _R(t) = R\exp (2\pi tM)$, where $M$ is a vector
in ${\cal I}_v$ with $|M|=1$, $0\le t\le 1/2$. Express $\tilde P$ through
variable $z$ also using $z^* = (2^v-2)^{-1} \{ -z
+ \sum_{s\in {\hat b}_v} s(zs^*) \} $.
Since $\lim_{|z|\to \infty } P(z) z^{-n-1}=1$, then due to Theorem 2.11
\cite{luoyst2}
$\lim_{R\to \infty } \int_{\gamma _R}(P{\tilde P})^{-1}(z)dz=$
$\int_{-R}^R(P{\tilde P})^{-1}(x)dx$ $=\int_{-R}^R |P(x)|^{-2}dx\ge 0.$
But $\lim_{R\to \infty } \int_{\gamma _R} (P{\tilde P})^{-1}(z)dz=
\lim_{R\to \infty }\pi R^{-2n-1}=0$.
On the other hand, $\int_{-R}^R|P(x)|^{-2}dx=0$ if and only if
$|P(x)|^{-2}=0$ for each $x\in \bf R$.
This contradicts our supposition, hence
there exists a zero $z_0\in {\cal A}_v$, that is, $P(z_0)=0$.
In the case $v=\infty $ use that $z=\lim_{v\to \infty }z_v$.

\thanks{
Address: Sergey V. Ludkovsky, Mathematical Department, TW-WISK,
Brussels University, V.U.B.,
Pleinlaan 2, Brussels 1050, Belgium. \\
{\underline {Acknowledgment}}. The author thanks the Flemish
Science Foundation for support through the Noncommutative Geometry
from Algebra to Physics project.}
\end{document}